\documentclass[conference,10pt,a4paper]{IEEEtran}
\usepackage{amssymb,amsmath,amsfonts,amstext,graphicx,amsthm,bm}
\title{Detection and estimation of spikes\\ in presence of noise and interference}
\author{\IEEEauthorblockN{Damien Passemier}
\IEEEauthorblockA{Department of Electronic\\ and Computer Engineering\\
Hong Kong University\\ of Science and Technology\\
Email: damien.passemier@gmail.com}
\and
\IEEEauthorblockN{Abla Kammoun}
\IEEEauthorblockA{CEMSE Division\\
King Abdullah University\\ of Science and Technology\\
Email: abla.kammoun@gmail.com}
\and
\IEEEauthorblockN{M\'erouane Debbah}
\IEEEauthorblockA{Alcatel-Lucent Chair\\ Sup\'elec\\
Gif-sur-Yvette, France\\
Email: merouane.debbah@supelec.fr}
\thanks{This research has been partly supported by the ERC Starting Grant 305123 MORE (Advanced Mathematical Tools for Complex Network Engineering).}}
\newtheorem{theorem}{Theorem}
\newtheorem{proposition}{Proposition}
\newtheorem{remark}{Remark}
\IEEEoverridecommandlockouts
\begin{document}
\maketitle
\begin{abstract}
	In many practical situations, the useful signal is contained in a low-dimensional subspace, drown in noise and interference. Many questions related to the estimation and detection of  the useful signal arise. Because of their particular structure, these issues are in connection to the problem that the mathematics community refers to as "spike detection and estimation". Previous works in this direction have been restricted to either determining the number of spikes  or estimating their values while knowing their multiplicities. This motivates our work which considers the joint estimation of the number of spikes and their corresponding orders, a problem which has not been yet investigated to the best of our knowledge.
%This makes our work interesting for both practical and mathematical perspectives.
%Beyond its importance for signal processing applications, this work is interesting from a mathematical perspective. 
	%is  proposes to consider the problem "spikes detection and estimation" from a different angle of view. 
	%Because of the particular structure of the signal, these issues  are in a close connection with the problem of "spike detection and estimation", as the mathematics community refers to it. 
\end{abstract}
\section{Introduction}

Detecting and estimating the components of a signal corrupted by additive Gaussion noise is a fundamental problem that arises in many signal and array processing applications.  Considering a large number of received samples, one can easy see that their covariance matrix  exhibit a different behaviour depending on the number of the components of the useful signal. In light of this consideration, first methods of signal detection like techniques using  the Roy Test \cite{roy} or  those using information theoretic criteria \cite{wax85} have been based on the eigenvalues of the empirical covariance matrix. Recently, the advances in the spectral analysis of large dimensional random matrices have engendered a new wave of interest for the scenario when the number of observations is of the same order of magnitude as the dimension of the received samples, while the number of signal components remain finite. Such a model is referred to as the spiked covariance model \cite{johnstone01}. It has allowed the emergence of new detection schemes based on the works of the extreme eigenvalues of large random Wishart matrices \cite{nadler10,nadler11,silverstein10}.  It is especially encountered in multi-sensor detection \cite{penna} and power estimation problems \cite{yao11}, which are at the heart of cognitive radio applications. This model has also found application in subspace estimation problems with a particular interest on the estimation of directions of arrival \cite{vallet}. 

From a mathematical perspective, the focus has been either to detect the presence of sources and estimate their numbers \cite{nadler09,passemier11} or to estimate their powers \cite{bai2012}.  The general case where the objective is to extract as much as possible information has not been addressed to the best of our knowledge. This motivates our work which proposes an easy way to jointly estimate the number of sources, their powers and their  multiplicities in the case where different sources are using the same power values.

\section{System model}
\label{sec:model}
Consider the $p$-dimensional observation vector ${\bf x}_i\in \mathbb{C}^p$ at time $i$: 
$$
{\bf x}_i=\sum_{k=1}^K \sqrt{\alpha_k}{\bf W}_k {\bf s}_{k,i}+ \sigma {\bf e}_i
$$
where
\begin{itemize}
	\item $({\bf W}_k)_{1\le k\le K}$ is an orthogonal family of rectangular unitary  $p\times m_k$ matrices (i.e, for $1\le k \le K, {\bf W}_k$ has orthogonal columns and ${\bf W}_k{\bf W}_j^{\mbox{\tiny H}}={1}_{k=j}{\bf I}_p$);
	\item $(\alpha_k)_{1\le k\le K}$ are $K$ positive distinct scalars such that $\alpha_1>\alpha_2>\cdots>\alpha_K$;
	\item $({\bf s}_{k,i})_{1\le k\le K} \in \mathbb{C}^{m_k\times 1}$  are independent random vectors with zero mean and variance $1$;
	\item ${\bf e}_i  \in \mathbb{C}^{p\times 1}$ is complex Gaussian distributed (i.e. ${\bf e}_i \sim\mathcal{C}\mathcal{N}({\bf 0},{\bf I})$) and represent the interference and noise signal;
    \item $\sigma^2$ is the strength of the noise.
	\end{itemize}
Therefore, we consider $K$ distinct powers $\alpha_k$, each of multiplicity $m_k$. Gathering $n$ observations ${\bf x}_1,\dots,{\bf x}_n$ into a $p\times n$ observation matrix ${\bf X}=\left[{\bf x}_1,\dots,{\bf x}_n\right]$, we obtain
	\begin{align*}
		{\bf X}&=\left[{\bf W}_1,\cdots,{\bf W}_K\right]\begin{bmatrix} \sqrt{\alpha_1}{\bf I}_{m_1}  &  & {\bf 0} \\
& \ddots &   \\
		{\bf 0}	   & &\sqrt{\alpha_K}{\bf I}_{m_K}
\end{bmatrix}\\
&\times
\begin{bmatrix}
{\bf s}_{1,1} & \cdots & {\bf s}_{1,n} \\
	\vdots & \ddots & \vdots \\
	{\bf s}_{K,1}& \cdots & {\bf s}_{K,n}
\end{bmatrix} +\sigma \left[{\bf e}_1,\cdots, {\bf e}_n\right].
\end{align*}
or equivalently
	\begin{equation}
	{\bf X}={\boldsymbol{\Sigma}}^{\frac{1}{2}}{\bf Y}
	\label{eq:model}
\end{equation}
	where ${\bf Y}$ is a matrix of independent entries with zero mean and variance $1$ and ${\bf \Sigma}$ is the theoretical covariance matrix of the observations given by:
	$$
	\boldsymbol{\Sigma}=\mathbf{U}\begin{bmatrix} 
		(\alpha_1+\sigma^2) {\bf I}_{m_1} & &\ &{\bf 0} \\
																					& \ddots & &\\
												& & (\alpha_K+\sigma^2) {\bf I}_{m_K}& \\
									{\bf 0}   & & & \sigma^{2}{\bf I}_{n-m}
	\end{bmatrix}\mathbf{U}^{\mbox{\tiny H}}
	$$
	where $m=\sum_{k=1}^K m_k$ and $\mathbf{U}$ is an orthogonal matrix. Note that $\boldsymbol{\Sigma}$ 
	has $K$ distinct eigenvalues with multiplicities $m_1,\dots,m_K$ and one eigenvalue equal to $\sigma^2$ with multiplicity $p-m$. 

	This model corresponds to the spiked covariance model \cite{johnstone01}: here we allow spikes with multiplicities greater than one. It can be encountered as shown in \cite{yao11} for power estimation purposes in cognitive radio networks. Another interesting application is met in the array processing field and in particular in the problem of the estimation of the angles of arrival. In this case,  the received signal matrix is given by \cite{vallet}:
	\begin{equation}
		{\bf X}={\bf A}({\bf \theta}){\bf P}^{\frac{1}{2}}{\bf S}+\sigma{\bf N}
	\label{eq:DOA}
\end{equation}
	where ${\bf A}=\left[{\bf a}(\theta_1),\dots, {\bf a}(\theta_m)\right]$, ${\bf a}(\theta_i)$ being the $p \times 1$ steering vector, ${\bf P}={\rm diag}\left({\alpha}{\bf I}_{m_1},\dots,\alpha_K{\bf I}_{m_K}\right)$, ${\bf S}$ is the $m\times n$ transmitted matrix of i.i.d Gaussian entries and ${\bf N}=\left[{\bf e}_1,\cdots, {\bf e}_n\right]$. Note that in this case, ${\bf A}$ can be considered as unitary, since ${\bf A}^{\mbox{\tiny H}}{\bf A}\rightarrow {\bf I}_m$ when $p\rightarrow \infty$.
	
Previous methods dealing with the estimation of  directions of arrivals has so far assumed a prior estimation of the number of sources \cite{vallet}. 
	Such information is obviously not always available in practice. This motivates our paper, which proposes a method to jointly estimate the number of sources as well as their multiplicities. %a mandatory piece of information on which rely several signal processing methods in practice.

	%Our objective in this paper is to estimate the multiplicity of each $\alpha_i$ as well as its value.
	\section{Estimation of spikes' values and  multiplicities}
	The estimation technique relies on  results about the asymptotic behavior of the covariance matrix. 
	As shown in the following proposition proven in \cite{bai08}, the asymptotic spectral properties of the covariance matrix depend on the eigenvalues $\alpha_1,\cdots,\alpha_K$ of the matrix $\boldsymbol{\Sigma}$. 
	
	%Prior to getting into the substance of the proposed method, we shall first recall the following results.
	\begin{proposition}
		Let ${\bf S}_n$ be the sample covariance matrix given by:
		$$
		{\bf S}_n=\frac{1}{n}\sum_{k=1}^n {\bf x}_i{\bf x}_i^{\mbox{\tiny H}}
		$$
		Denote by $\widehat{\lambda}_{n,1}>\widehat{\lambda}_{n,2} >\cdots > \widehat{\lambda}_{n,p} $ the $p$ eigenvalues of ${\bf S}_n$  arranged in decreasing order. Let $s_i=\sum_{k=1}^i m_k$ and $J_k$ the index set $J_k=\left\{s_k+1,\dots,s_k+m_k\right\}$, $k\in\left\{1,\dots,K\right\}$.
		
		Let $\phi(x)=x+\sigma^2+\gamma\sigma^2\left (1+\frac{\sigma^2}{x}\right )$ for $x\neq 0$ and assume that $\gamma_n=\frac{p}{n}\to \gamma$. 
		Then, if $\phi'(\alpha_k)>0$   (i.e. $\alpha_k> \sigma^2\sqrt{\gamma}$) for any $k\in\left\{1,\dots,K\right\}$, we have almost surely
		$$
		\widehat{\lambda}_{n,j}\to \phi(\alpha_k),  \ \ \forall j\in J_k
		$$
	\end{proposition}
	\begin{remark}
 Under the condition $\phi'(\alpha_k)>0$ for all $k\in\left\{1,\dots,K\right\}$, the empirical distribution of the spectrum is composed of $K+1$ connected intervals: a bulk corresponding to the Mar{\v c}henko-Pastur law \cite{baik06} followed by $K$ spikes. To illustrate this, we represent in Figure \ref{fig:histogram}, the empirical histogram of the eigenvalues of the empirical covariance matrix when $K=3$, $(\alpha_1,\alpha_2,\alpha_3)=(7,5,3)$, $(n,p)=(4000,2000)$ and $\sigma^2=1$. 
		\begin{figure}
			\begin{center}
				\includegraphics[scale=0.4]{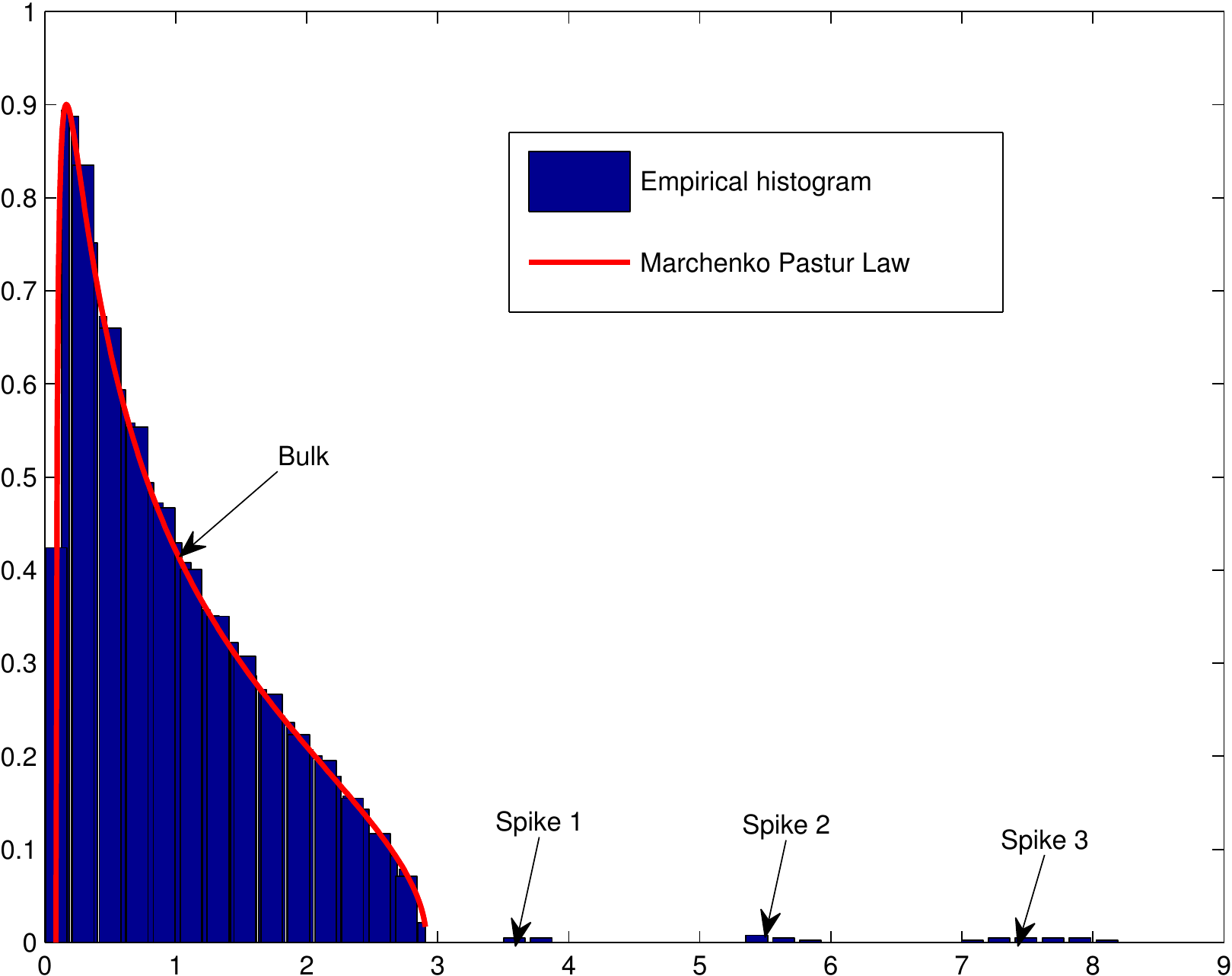}
			\end{center}
			\caption{Histogram of the eigenvalues of the empirical covariance matrix.}
			\label{fig:histogram}
		\end{figure}

	\end{remark}
	Figure \ref{fig:histogram} provides us insights about an  intuitive approach to estimate the multiplicities of spikes and their values given their number $K$. Actually, one needs to rearrange the eigenvalues and then detect the largest gaps that correspond to a switch from one connected interval to the next one. 

This leads us to distinguish two cases whether $K$ is either known or not. We will consider these cases in turn in the following.
	%From fig. \ref{fig:histogram}, an intuitive way to estimate the multiplicities of spikes given their number $K$, consisting in rearranging the empirical eigenvalues in order to detect the $K$ largest gaps that correspond to a switch from a spike to another one. 
%	We will distinguish two cases whether $K$ is known or not.
	\subsection{$K$ is known}
	\label{sec:known}
In this case, we propose to estimate the eigenvalues by considering the differences between consecutive eigenvalues:
$$
\delta_{n,j}=\widehat{\lambda}_{n,j}-\widehat{\lambda}_{n,j+1},  \ \ j\geq 1.
$$
Indeed, the results quoted above imply that a.s. $\delta_{n,j} \rightarrow 0$, for $j \notin \{s_i\text{, } i=1,\dots,K\}$ whereas for $j \in \{s_i\text{, } i=1,\dots,K\}$, $\delta_{n,j}$ tends to a positive limit given by $\phi(\alpha_j)-\phi(\alpha_{j+1})$. Thus it becomes possible to estimate the multiplicities from index-numbers $j$ where $\delta_{n,j}$ is large. If $K$ is known, we will take the indices corresponding to the $K$ larger differences $\delta_{n,i}$. Denote by $i_1,\dots,i_p$ the indices of the differences $\delta_{n,i}$ such that $\delta_{n,i_1} \ge \dots \ge \delta_{n,i_p}$. Then, the estimator ($\hat{m}_1,\dots,\hat{m}_K$) of the multiplicities ($m_1,\dots,m_K$) is defined by
\[ \left \{
\begin{array}{@{}l@{}c@{}l}
\hat{m}_1&=&\text{min~} \{ i_k\text{, }k \in \{1,\dots,K\}\}\vspace{0.3cm}\\
\hat{m}_2&=&\text{min~} \{ i_k\text{, }k \in \{1,\dots,K\}\backslash \{\hat{m}_1\}\}-\hat{m}_1\\
\hat{m}_j&=&\text{min~} \{ i_k\text{, }k \in \{1,\dots,K\}\backslash \{\hat{m}_1,\dots,\hat{m}_{j-1}\}\}-\displaystyle\sum_{i=1}^{j-1}\hat{m}_i\\
\hat{m}_K&=&\text{max~} \{ i_k\text{, }k \in \{1,\dots,K\}\}-\displaystyle\sum_{i=1}^{K-1}\hat{m}_i
\end{array}\right.\]
%Notice that $\text{max~} \{ i_k\text{, }k \in \{1,\dots,K\}\}$ is also an estimator of the number of spikes.
The proposed consistent estimator of the number of the spikes is therefore given by the following theorem, for which a proof is omitted because of lack of space:
\begin{theorem} \label{consistency}
	Let  $(\mathsf{x}_i)_{1 \le i \le n}$ be $n$ i.i.d copies. of $\mathsf{x}$ which follows the model \eqref{eq:model}. Suppose that the population covariance matrix $\Sigma$ has $K$ non null eigenvalues $(\alpha_i+\sigma^2)_{1 \le i \le K}$ such that $\alpha_1>\cdots>\alpha_{K}>\sigma^2\sqrt{\gamma}$ with respective multiplicity $(m_k)_{1\le k \le K}$ ($m_1+\dots+m_K=m$), and $p-m$  eigenvalues equal to $\sigma^2$. Assume that $p/n \rightarrow \gamma >0$ when $n \rightarrow \infty$. Then the estimator ($\hat{m}_1,\dots,\hat{m}_K$) is strongly consistent, i.e $(\hat{m}_1,\dots,\hat{m}_K)  \rightarrow (m_1,\dots,m_K)$ almost surely when $n \rightarrow \infty$.
\end{theorem}
\subsection{$K$ is not known}
As Figure \ref{fig:histogram} shows, eigenvalues outside the bulk are organized into $K$ clusters, where within each cluster, all eigenvalues converge  to the same value in the asymptotic regime $p,n\to+\infty$ such that $p/n\to\gamma$. If $K$ is not estimated correctly, applying the previous method, will lead to either gathering two close clusters ($K$ is under-estimated) or to subdividing the clusters corresponding to the highest spikes ($K$ is over-estimated). Clearly, the second order results within each cluster seems to bring useful information which allows to discard these cases. In particular, in the sequel, we will rely on the following proposition which is a by-product of Proposition 3.2 in \cite{bai08}:
\begin{proposition}
	Assume that the settings of Theorem \ref{consistency} holds.
	Let $g_k=\sum_{j=s_{k-1}+1}^{s_k} \widehat{\lambda}_{n,j}$, the sum of the  eigenvalues corresponding to the $k$-th cluster. Then, when $n \rightarrow \infty$ such that $p/n \rightarrow \gamma >0$, $g_k$ verify
	$$
	\sqrt{n}\left(g_k-m_k\phi(\alpha_k)\right)\xrightarrow{\mathcal{L}}\mathcal{N}(0,2m_kv_k^2)
	$$
	where $v_k^2=2\alpha_k'^2((\alpha_k'-1)^2-\gamma)/(\alpha_k'-1)^2$, $\alpha_k'=\alpha_k/\sigma^2+1$ and $\mathcal{L}$ denotes the convergence in distribution.
	\label{th:clt}
\end{proposition}
Theorem \ref{th:clt} establishes that the sum of the eigenvalues within the $k$-th cluster behaves as a Gaussian random variable with  mean and variance depending on the unknown value $\alpha_k$. One way to remove the uncertainty in the unknowns $\alpha_k$ is to assume that they are random with a priori known distribution $\pi\left(\alpha_1,\dots,\alpha_K|K\right)$. A possible case would correspond to the situation where they are uniformly distributed over a finite discrete set\footnote{A discrete distribution for powers has been considered in \cite{HAC04}.}.

Since the clusters are asymptotically independent \cite{couillet13}, the likelihood function (distribution of  ${\bf g}=\left[g_1,\cdots,g_K\right]$ under the underlying parameters $\alpha_1,\cdots,\alpha_K,m_1,\cdots,m_K, K$)  is given by:
$$
f({\bf g}|\alpha_1,\dots,\alpha_K,K)=\prod_{k=1}^K\frac{1}{\sqrt{2\pi v_k^2}}e^{-\frac{1}{2v_k^2}\left(g_k-m_k\phi(\alpha_k)\right)^2}
$$
where the multiplicities $m_1,\dots,m_K$ can be estimated in a consistent way given the number of classes $K$ as it has been shown in  section \ref{sec:known}.
Hence, the maximum likelihood function $f({\bf g}|K)$ is given by:
\begin{equation}
f({\bf g}|K)=\mathbb{E}\left[f({\bf g}|\alpha_1,\dots,\alpha_K,K)\right]
\label{eq:mle}
\end{equation}
where the expectation is taken over the a priori distribution $\pi(\alpha_1,\cdots,\alpha_K|K)$.
%As it was shown earlier, given $K$, the multiplicities $m_1,\cdots,m_K$ can be estimated in a consistent way.
The maximum likelihood estimator $\widehat{K}$ is thus given by
$$
\widehat{K}=\max_{1\leq k \leq K_{\rm max}} \mathbb{E}\left[f({\bf g}|\alpha_1,\dots,\alpha_K,K) \right]\text{,}
$$
where $K_{\rm max}$ is a known upper bound of $K$. Once $K$ is estimated, the multiplicities can be retrieved by using the method in Section \ref{sec:known}.

To sum up,  when $K$ is unknown, the estimation of  the unknown parameters using the a priori $\pi$ consists in the following steps :
\begin{enumerate}
\item Compute the consecutive differences of the ordered eigenvalues of the sample covariance matrix ${\bf S}_n$ given by $\delta_{n,j}=\widehat{\lambda}_{n,j}-\widehat{\lambda}_{n,j+1}$;
%\item Compute the maximum likelihood function given in  \eqref{eq:mle}.
\item	For each $k$ ranging from one to $K_{\rm max}$, calculate the corresponding estimator $(\hat{m}_1^{(k)},\dots,\hat{m}_k^{(k)})$ of the multiplicities using Theorem \ref{consistency}, and compute the maximum likelihood function \eqref{eq:mle}. 
%\item %For each $k$, calculate all the possible values of $(\alpha_1,\dots,\alpha_k)$ when taking $k$ spikes among the possible set of a priori values $E$. For each possibility, calculate the value of the likelihood $f$ using the corresponding multiplicities found step 2. Add these likelihoods together and divide by the number of possibilities (i.e. $|E|$ choose $k$) to obtain $R(K)$;
\item Select  $K$ such that it maximizes the maximum likelihood function.
\end{enumerate}

\section{Numerical experiments}
We consider in our simulations the model described by \eqref{eq:DOA} given in section \ref{sec:model} with ${\bf A(\theta)}=p^{-1/2}\left [\exp\left (-i v \sin(\theta) \pi \right )\right ]_{v=0}^{p-1}$, where $\theta$ is chosen uniformly on $[0,2\pi)$. We assume that the set of the a priori spikes is $E=\left\{1,3,5,7\right\}$ and that the values $\alpha_1,\dots,\alpha_K$ are uniformly distributed over this set. 
%To assess the  of our estimator $\hat{K}$, we perform some simulation experiments. We consider the model \eqref{eq:DOA} given Section \ref{sec:model}, with ${\bf A(\theta)}=p^{-1/2}\left [\exp\left (-i v \sin(\theta) \pi \right )\right ]_{v=0}^{p-1}$. We assume that the set of the a priori spikes is $E=\left \{ 1, 3, 5, 7 \right \}$ and we fix the unknown $K$. We assume that $m$ is known and take $K_{\rm max}=m$.

In the sequel, we will display the empirical probability $\mathbb{P}(\hat{K}=K)$ calculated over $500$ independent realizations. For each iteration, we choose the ``true'' values of spikes uniformly in the set $E$, but with the same fixed proportion $m_i/m$, $i=1,\dots,K$.

We consider two different experiments: in the first one,  we study the performance of our method for different level of noise variances  whereas for the second one, we consider the impact of the number of spikes $m$ for a fixed noise variance.

\subsection{Performance of the proposed method with respect to the variance of the noise}
In this experiment, we consider the detection of the number of $K=3$ different clusters of $500\times 1$ ( $p=500$ ) signals from $n=1000$ samples. We assume that the unknown multiplicities are $m_1=1$, $m_2=4$, $m_3=2$. Since the minimum value of the spike is assumed to be 1, $\sigma^2$ has to be lower than $1/\sqrt{c}=1.4142$ in order to keep a gap between $\hat{\lambda}_m$ and $\hat{\lambda}_{m+1}$ (see Theorem 1). The noise variance is expressed in dB $10\log_{10}(\sigma^2)$.  Table \ref{tab:SNR} illustrates the obtained results : 
%		\begin{figure}
%			\begin{center}
%				\includegraphics[scale=0.8]{SNR.pdf}
%			\end{center}
%\vspace{-0.8cm}
%			\caption{Empirical probability of $\mathbb{P}(\hat{K}=K)$ as a function of the SNR.}
%			\label{fig:SNR}
%		\end{figure}

{\scriptsize
\begin{table}[!ht]
\caption{{  Empirical probability of $\mathbb{P}(\hat{K}=K)$ as a function of the $\sigma^2$.}}
{\scriptsize\begin{tabular}{|c|c@{~~}c@{~~}c@{~~}c@{~~}c@{~~}c@{~~}c@{~~}c@{~~}c@{~~}c@{~~}|}
\hline
$\sigma^2$(dB) & -50 & -40 & -30 & -20 & -10 & -6.99 & -5.223 & -3.98 & -3.01 & -2.22 \\
\hline
$\mathbb{P}(\hat{K}=K)$ & 0.992 & 0.978 & 0.988 & 0.986 & 0.984  & 0.978 & 0.978 & 0.980 & 0.964 & 0.974 \\
\hline
\hline
SNR (dB) & -1.55 & -0.97 & -0.46 & 0 & 0.41 & 0.80 & 0.97 & 1.14 & 1.30 & 1.46\\
\hline
$\mathbb{P}(\hat{K}=K)$ & 0.972 & 0.954 & 0.960 & 0.968 & 0.942 & 0.926 & 0.896 & 0.850 & 0.694 & 0.476  \\
\hline
\end{tabular}}
\label{tab:SNR}

\end{table}}
%\vspace{-0.3cm}

Our estimator performs well, especially for low noise variances. When $\sigma^2$ is getting close to the threshold 1.41 (i.e. 1.50 dB), the estimator becomes less accurate, which was expected since $\hat{\lambda}_m$ is very close to the bulk.

\subsection{Influence of the number of spikes $m$}
\begin{figure}
			\begin{center}
				\includegraphics[scale=0.8]{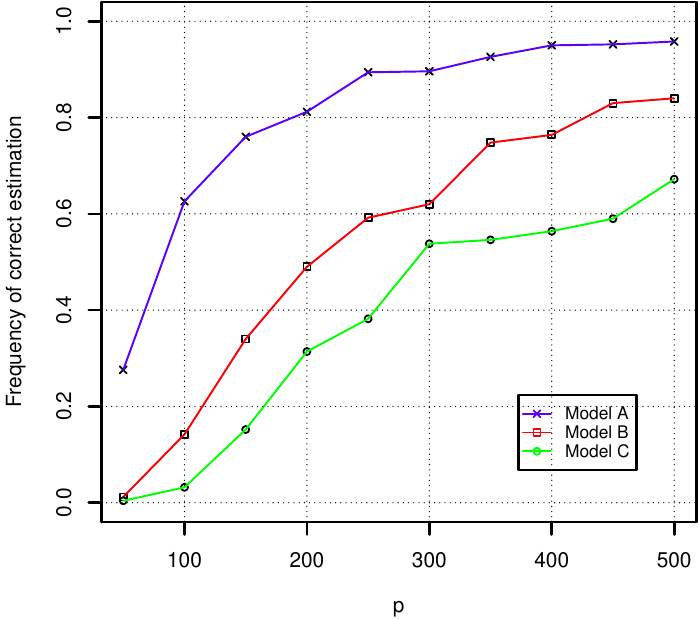}
			\end{center}
\vspace{-0.7cm}
			\caption{Empirical probability of $\mathbb{P}(\hat{K}=K)$ as a function of $(p,n)$, for Models A, B and C.}
			\label{fig:pn}
		\end{figure}

%Here, we fix $\sigma^2=1$. We consider three models, all with $K=3$:
We study in this experiment the impact of the number of spikes in the performance of the proposed estimation method. Similarly to the previous simulation settings, we set $K=3$ and $\gamma=p/n=0.5$. We consider the following three models:
\begin{itemize}
\item Model A: $m=4$, with $m_1=1$, $m_2=2$, $m_3=1$;
\item Model B: $m=8$, with $m_1=2$, $m_2=4$, $m_3=2$;
\item Model C: $m=12$, with $m_1=3$, $m_2=6$, $m_3=3$;
\end{itemize}
Figure \ref{fig:pn} displays the frequency of correct estimation for these three models with respect to $p$. Note that these models keep the $p/n$ and $m_i/m$ fixed except $m/n$ which is different. In that way, only the impact of the variation of the number of spikes is visualized. 

%$c=0.5$ is kept fixed and $p$ (or equivalently $n$) is varying. The results are displayed Fig. \ref{fig:pn}.

As expected, our estimator performs better in Model A than in Model B and C. In both cases, we observe the asymptotic consistency, but the convergence is slower for Model C.
\begin{remark}
Once $K$ was correctly estimated, we have noticed by simulations that the multiplicities are correctly estimated. This is in accordance with our Theorem \ref{consistency}.
\end{remark}

\section{Conclusion}
%In this paper we have considered the problem of estimating the multiplicities of the spikes in a spiked population model, which has not been considered yet to the best of our knowledge. When the number of different spikes $K$ is known, we have provided a consistent estimator of the multiplicities. When $K$ is unknown, we have given an estimator of $K$ based on a likelihood estimation. Our simulation experiments have illustrated the accuracy of the estimation in different settings.
The problem of signal detection appears naturally in many signal processing applications. Previous works used to deal with this problem partially by assuming extra knowledge about the number of spikes or their corresponding orders. This work is therefore an attempt to consider the general problem where the objective is to estimate all the unknown parameters. In particular, we show that when the number of different spikes is known, their multiplicities  can be estimated consistently. In light of this consideration, we propose a Bayesian estimation method which jointly infer the number of spikes and their multiplicities. The  experiments that we carried out support the performance of the proposed technique.

%Experiments support the performance of the proposed method. 

%We assume thus that the receiver has a prior distribution of the $\alpha_k, k\in\left\{1,\cdots,K\right\}.$

%More formally, assume that $\alpha_1,\cdots,\alpha_K$ are random with a prior probability distribution $\pi(\alpha_1,\cdots,\alpha_K|K)$, known at the receiver. 

%
%As far as the complex case is considered, the joint eigenvalue distribution of the eigenvalue within each cluster is not 

\bibliographystyle{IEEEtrans}
\bibliography{mybib}

% Generated by IEEEtranS.bst, version: 1.13 (2008/09/30)
\begin{thebibliography}{10}
\providecommand{\url}[1]{#1}
\csname url@samestyle\endcsname
\providecommand{\newblock}{\relax}
\providecommand{\bibinfo}[2]{#2}
\providecommand{\BIBentrySTDinterwordspacing}{\spaceskip=0pt\relax}
\providecommand{\BIBentryALTinterwordstretchfactor}{4}
\providecommand{\BIBentryALTinterwordspacing}{\spaceskip=\fontdimen2\font plus
\BIBentryALTinterwordstretchfactor\fontdimen3\font minus
  \fontdimen4\font\relax}
\providecommand{\BIBforeignlanguage}[2]{{%
\expandafter\ifx\csname l@#1\endcsname\relax
\typeout{** WARNING: IEEEtranS.bst: No hyphenation pattern has been}%
\typeout{** loaded for the language `#1'. Using the pattern for}%
\typeout{** the default language instead.}%
\else
\language=\csname l@#1\endcsname
\fi
#2}}
\providecommand{\BIBdecl}{\relax}
\BIBdecl

\bibitem{bai2012}
Z.~Bai and X.~Ding, ``{Estimation of Spiked Eigenvalues in Spiked Models},''
  \emph{Random Matrices: Theory and Applications}, vol.~1, no.~2, p. 1150011,
  2012.

\bibitem{bai08}
Z.~D. Bai and J.-F. Yao, ``{Central Limit Theorems for Eigenvalues in a Spiked
  Population Model},'' \emph{{Annales de l'Institut Henri Poincar{\'e}}},
  vol.~44, no.~3, pp. 447--474, 2008.

\bibitem{baik06}
J.~Baik and J.~W. Silverstein, ``{Eigenvalues of Large Sample Covariance
  Matrices of Spiked Population Models},'' \emph{Journal Of Multivariate
  Analysis}, vol.~97, no.~6, pp. 1382--1408, 2006.

\bibitem{HAC04}
J.-M. Chaufray, W.~Hachem, and P.~Loubaton, ``{Asymptotic Analysis of Optimum
  and Sub-Optimum CDMA Sownlink MMSE Receivers},'' \emph{IEEE Transactions on
  Information Theory}, vol.~50, no.~11, pp. 2620--2638, 2004.

\bibitem{vallet}
W.~Hachem, P.~Loubaton, X.~Mestre, J.~Najim, and P.~Vallet, ``{A Subspace
  Estimator for Fixed Rank Perturbations of Large Random Matrices},''
  \emph{Journal of Multivariate Analysis}, vol. 114, pp. 427--447, 2013.

\bibitem{johnstone01}
I.~M. Johnstone, ``{On the Distribution of the Largest Eigenvalue in Principal
  Component Analysis},'' \emph{{Annals of Statistics}}, vol.~29, no.~2, pp.
  295--327, 2001.

\bibitem{nadler09}
S.~Kritchman and B.~Nadler, ``{Non-Parametric Detection of the Number of
  Signals: Hypothesis Testing and Random Matrix Theory},'' \emph{IEEE
  Transactions on Signal Processing}, vol.~57, no.~10, pp. 3930--3941, 2009.

\bibitem{silverstein10}
R.~R. Nadakuditi and J.~W. Silverstein, ``{Fundamental Limit of Sample
  Generalized Eigenvalue Based Detection of Signals in Noise Using Relatively
  few Signal-bearing and noise-only samples},'' \emph{{IEEE Selected Topics in
  Signal Processing}}, vol.~4, no.~3, pp. 468--480, 2010.

\bibitem{nadler10}
B.~Nadler, ``{Nonparametric Detection of Signals by Information Theoretic
  Criteria: Performance Analysis and an Improved Estimator},'' \emph{{IEEE
  Transactions on Signal Processing}}, vol.~58, no.~5, pp. 2746--2756, 2010.

\bibitem{nadler11}
------, ``{Detection Performance of Roy's Largest Root Test When the Noise
  Covariance Matrix is Arbitrary},'' in \emph{SSP}, Nice, 2011.

\bibitem{passemier11}
D.~Passemier and J.~F. Yao, ``{Estimation of the Number of Spikes, Possibly
  Equal, in the High- Dimensional Case},'' \emph{Journal Of Multivariate
  Analysis}, vol. 127, pp. 173--183, 2014.

\bibitem{penna}
F.~Penna, ``{Statistical Methods for Cooperative and Distributed Inference in
  Wireless Networks},'' Ph.D. dissertation, Politecnico di Torino, Torino,
  Italy, 2012.

\bibitem{couillet13}
{R. Couillet and W. Hachem}, ``{Fluctuations of Spiked Random Matrix Models and
  Failure Diagnosis in Sensor Networks},'' \emph{IEEE Transactions on
  Information Theory}, vol.~59, no.~1, pp. 509--525, 2013.

\bibitem{roy}
S.~N. Roy, \emph{{Some Aspects of Mutlivariate Analysis}}.\hskip 1em plus 0.5em
  minus 0.4em\relax New-York: Wiley, 1957.

\bibitem{wax85}
M.~Wax and T.~Kailath, ``{Detection of Signals by Information Theoretic
  Criteria},'' \emph{{IEEE Transactions on Acoustics, Speech and Signal
  Processing}}, vol.~33, no.~2, pp. 387--392, 1985.

\bibitem{yao11}
J.~Yao, R.~Couillet, J.~Najim, E.~Moulines, and M.~Debbah, ``{CLT for
  Eigen-Inference Methods for Cognitive Radio},'' in \emph{ICASSP}, Prague,
  2011.

\end{thebibliography}
%\bibliography{biblio}

\end{document}